# SPACELIKE $B_2$-SLANT HELIX IN MINKOWSKI 4-SPACE $E_1^4$

**Mehmet ÖNDER, Hüseyin KOCAYİĞİT, Mustafa KAZAZ**
*Celal Bayar University, Faculty of Science and Art, Department of Mathematics, 45047 Manisa, Turkey.* e-mail: mehmet.onder@bayar.edu.tr

## ABSTRACT

In this paper, we give the characterizations of spacelike $B_2$-slant helix by means of curvatures of the spacelike curve in Minkowski 4-space. Furthermore, we give the integral characterization of the spacelike $B_2$-slant helix.

## 1. INTRODUCTION

A curve of constant slope or general helix in Euclidean 3-space $E^3$ is defined by the property that the tangent makes a constant angle with a fixed straight line (the axis of the general helix)[2]. A classical result stated by M. A. Lancret in 1802 and first proved by B. de Saint Venant in 1845 (see [12] for details) is: *A necessary and sufficient condition that a curve be a general helix is that the ratio of the first curvature to the second curvature be constant i.e., $k_1/k_2$ is constant along the curve, where $k_1$ and $k_2$ denote the first and second curvatures of the curve, respectively.* Analogue to that A. Magden has given a characterization for a curve $x(s)$ to be a helix in Euclidean 4-space $E^4$. He characterizes a helix iff the function

$$\frac{k_1^2}{k_2^2} + \left[\frac{1}{k_3}\frac{d}{ds}\left(\frac{k_1}{k_2}\right)\right]^2$$

is constant where $k_1, k_2$ and $k_3$ are first, second and third curvatures of Euclidean curve $x(s)$, respectively, and they are nowhere zero[7]. Corresponding characterizations of timelike helices in Minkowski 4-space $E_1^4$ were given by Kocayigit and Onder [5]. Latterly Camci and et al. have given some characterizations for a non-degenerate curve a to be a generalized helix by using its harmonic curvatures [2].

Recently, Izumiya and Takeuchi have introduced the concept of slant helix by saying that the normal lines of the curve make a constant angle with a fixed direction and they have given a characterization of slant helix in Euclidean 3-space $E^3$ by the fact that the function

$$\frac{k_1^2}{(k_1^2+k_2^2)^{3/2}}\left(\frac{k_2}{k_1}\right)'$$

is constant [4]. After them, Kula and Yayli investigated spherical images, the tangent indicatrix and the binormal indicatrix of a slant helix and they obtained that the spherical images are spherical helices [6]. Analogue to the definition of slant helix, Onder and et al. have defined $B_2$-slant helix in Euclidean 4-space $E^4$ by saying that the second binormal vector of the curve make a constant angle with e fixed direction and they have given some characterizations of $B_2$-slant helix in Euclidean 4-space $E^4$ [8].

In this paper, we consider spacelike $B_2$-slant helix in Minkowski 4-space $E_1^4$ and we give some characterizations and also the integral characterization of spacelike $B_2$-slant helix.

## 2. PRELIMINARIES

Minkowski space-time $E_1^4$ is an Euclidean space $E^4$ provided with the standard flat metric given by

$$\langle,\rangle = -dx_1^2 + dx_2^2 + dx_3^2 + dx_4^2$$

where $(x_1, x_2, x_3, x_4)$ is a rectangular coordinate system in $E_1^4$.

Since $\langle,\rangle$ is an indefinite metric, recall that a vector $v \in E_1^4$ can have one of three causal characters; it can be spacelike if $\langle v,v \rangle > 0$ or $v = 0$, timelike if $\langle v,v \rangle < 0$ and null(lightlike) if $\langle v,v \rangle = 0$ and $v \neq 0$. Similarly, an arbitrary curve $x(s)$ in $E_1^4$ can locally be spacelike, timelike or null (lightlike), if all of its velocity vectors $x'(s)$ are respectively spacelike, timelike or null (lightlike). Also recall that the pseudo-norm of an arbitrary vector $v \in E_1^4$ is given by $\|v\| = \sqrt{|\langle v,v \rangle|}$. Therefore $v$ is a unit vector if $\langle v,v \rangle = \pm 1$. The velocity of the curve $x(s)$ is given by $\|x'(s)\|$. Next, vectors $v$, $w$ in $E_1^4$ are said to be orthogonal if $\langle v,w \rangle = 0$. We say that a timelike vector is *future pointing* or *past pointing* if the first compound of the vector is positive or negative, respectively.

Denote by $\{T(s), N(s), B_1(s), B_2(s)\}$ the moving Frenet frame along the curve $x(s)$ in the space $E_1^4$. Then $T$, $N$, $B_1$, $B_2$ are the tangent, the principal normal, the first binormal and the second binormal fields, respectively. A timelike (resp. spacelike) curve $x(s)$ is said to be parameterized by a pseudo-arclength parameter $s$, i.e. $\langle x'(s), x'(s) \rangle = -1$ (resp. $\langle x'(s), x'(s) \rangle = 1$ ).

Let $x(s)$ be a spacelike curve in Minkowski space-time $E_1^4$, parameterized by arclength function of $s$. Then for the curve $x(s)$ the following Frenet equation is given as follows

$$\begin{bmatrix} T' \\ N' \\ B_1' \\ B_2' \end{bmatrix} = \begin{bmatrix} 0 & k_1 & 0 & 0 \\ -\varepsilon_1 k_1 & 0 & k_2 & 0 \\ 0 & \varepsilon_2 k_2 & 0 & k_3 \\ 0 & 0 & \varepsilon_1 k_3 & 0 \end{bmatrix} \begin{bmatrix} T \\ N \\ B_1 \\ B_2 \end{bmatrix}, \quad (1)$$

where

$$\langle T,T \rangle = 1, \ \langle N,N \rangle = \varepsilon_1, \ \langle B_2, B_2 \rangle = \varepsilon_2, \ \langle B_1, B_1 \rangle = -\varepsilon_1 \varepsilon_2, \ \varepsilon_1 = \pm 1, \ \varepsilon_2 = \pm 1$$

and recall that the functions $k_1 = k_1(s)$, $k_2 = k_2(s)$ and $k_3 = k_3(s)$ are called the first, the second and the third curvature of the spacelike curve $x(s)$, respectively and we will assume throughout this work that all the three curvatures satisfy $k_i(s) \neq 0$, $1 \leq i \leq 3$. Here the signs of $\varepsilon_1$ and $\varepsilon_2$ are changed by a rule. The signature rule between $\varepsilon_1$ and $\varepsilon_2$ can be given as follows

|     | $\varepsilon_1$ |      | $\varepsilon_2$ |
|-----|------|------|-----------|
| if  | +1   | then | +1 or -1  |
|     | -1   |      | +1        |

or

|     | $\varepsilon_2$ |      | $\varepsilon_1$ |
|-----|------|------|-----------|
| if  | +1   | then | +1 or -1  |
|     | -1   |      | +1        |

For the obvious forms of the Frenet equations in (1) we refer to the reader to see ref. [14].

## 3. SPACELIKE $B_2$-SLANT HELICES IN MINKOWSKI 4-SPACE

In this section, we give the definition and the characterizations of spacelike $B_2$-slant helix.

Let $x : I \subset IR \to E_1^4$ be a unit speed spacelike curve with nonzero curvatures $k_1(s)$, $k_2(s)$ and $k_3(s)$ and let $\{T, N, B_1, B_2\}$ denotes the Frenet frame of the curve $x(s)$. We call $x(s)$ as spacelike $B_2$-slant helix if its second binormal unit vector $B_2$ makes a constant angle with a fixed direction in a unit vector $U$; that is

$$\langle B_2, U \rangle = constant \tag{2}$$

along the curve. By differentiation (2) with respect to $s$ and using the Frenet formulae (1) we have

$$\langle \varepsilon_1 k_3 B_1, U \rangle = 0.$$

Therefore $U$ is in the subspace $Sp\{T, N, B_2\}$ and can be written as follows

$$U = a_1(s)T(s) + a_2(s)N(s) + a_3(s)B_2(s), \tag{3}$$

where

$$a_1 = \langle U, T \rangle, \quad \varepsilon_1 a_2 = \langle U, N \rangle, \quad \varepsilon_2 a_3 = \langle U, B_2 \rangle = constant.$$

Since $U$ is unit, we have

$$a_1^2 + \varepsilon_1 a_2^2 + \varepsilon_2 a_3^2 = M. \tag{4}$$

Here $M$ is +1, -1 or 0 depending if $U$ is spacelike, timelike or lightlike, respectively. The differentiation of (3) gives

$$\left(\frac{da_1}{ds} - \varepsilon_1 a_2 k_1\right)T + \left(\frac{da_2}{ds} + a_1 k_1\right)N + (a_2 k_2 + \varepsilon_1 a_3 k_3)B_1 = 0,$$

and from this equation we get

$$a_2 = -\varepsilon_1 \frac{k_3}{k_2} a_3 = \varepsilon_1 \frac{1}{k_1} \frac{da_1}{ds}, \quad \frac{da_2}{ds} = -a_1 k_1. \tag{5}$$

Since

$$\frac{da_2}{ds} = -a_1 k_1 \quad \text{and} \quad \frac{da_2}{ds} = -\varepsilon_1 \frac{k_1'}{k_1^2} \frac{da_1}{ds} + \frac{\varepsilon_1}{k_1} \frac{d^2 a_1}{ds^2},$$

we find the second order linear differential equation in $a_1$ given by

$$\varepsilon_1 \frac{d^2 a_1}{ds^2} - \varepsilon_1 \frac{k_1'}{k_1} \frac{da_1}{ds} + a_1 k_1^2 = 0. \tag{6}$$

If we change variables in the above equation as $t = \int_0^s k_1(s) ds$ then we get

$$\frac{d^2 a_1}{dt^2} + \varepsilon_1 a_1 = 0.$$

This equation has two solutions: If $\varepsilon_1 = +1$ then the solution is

$$a_1 = A \cos \int_0^s k_1(s) ds + B \sin \int_0^s k_1(s) ds, \tag{7}$$

and if $\varepsilon_1 = -1$ then the solution is

$$a_1 = A \cosh \int_0^s k_1(s) ds + B \sinh \int_0^s k_1(s) ds. \tag{8}$$

where $A$ and $B$ are constant.

Assume that $\varepsilon_1 = -1$ and consider the solution (8). From (5) and (8) we have

$$a_2 = \frac{k_3}{k_2} a_3 = -A \sinh \int_0^s k_1(s) ds - B \cosh \int_0^s k_1(s) ds,$$

$$a_1 = -\frac{1}{k_1} \left( \frac{k_3}{k_2} \right)' a_3 = A \cosh \int_0^s k_1(s) ds + B \sinh \int_0^s k_1(s) ds.$$

From these equations it follows that

$$A = a_3 \left( \frac{k_3}{k_2} \sinh \int_0^s k_1(s) ds - \frac{1}{k_1} \left( \frac{k_3}{k_2} \right)' \cosh \int_0^s k_1(s) ds \right), \tag{9}$$

$$B = a_3 \left( -\frac{k_3}{k_2} \cosh \int_0^s k_1(s) ds + \frac{1}{k_1} \left( \frac{k_3}{k_2} \right)' \sinh \int_0^s k_1(s) ds \right). \tag{10}$$

Hence, using (9) and (10) we get

$$B^2 - A^2 = \left[ \left( \frac{k_3}{k_2} \right)^2 - \frac{1}{k_1^2} \left( \left( \frac{k_3}{k_2} \right)' \right)^2 \right] a_3^2 = constant,$$

so that

$$\left( \frac{k_3}{k_2} \right)^2 - \frac{1}{k_1^2} \left( \left( \frac{k_3}{k_2} \right)' \right)^2 = constant := m. \tag{11}$$

From (4), (9), (10) and (11) we have

$$B^2 - A^2 = a_3^2 m = M - 1.$$

Thus, the sign of the constant $m$ agrees with the sign of $B^2 - A^2$. So, if $U$ is timelike or lightlike then $m$ is negative. If $U$ is spacelike then $m = 0$. Then we can give the following corollary.

**Corollary 3.1.** *Let $x(s)$ be a spacelike $B_2$-slant helix with timelike principal normal $N$ in Minkowski 4-space $E_1^4$ and $U$ be a unit constant vector which makes a constant angle with the second binormal $B_2$. Then the vector $U$ is spacelike if and only if there exist a constant $K$ such that*

$$\frac{k_3}{k_2}(s) = K \exp\left(\int_0^s k_1(t)dt\right).$$

When $\varepsilon_1 = +1$, by using (7) with similar calculations as above we get that the spacelike curve $x(s)$ is a spacelike $B_2$-slant helix if and only if

$$\left(\frac{k_3}{k_2}\right)^2 + \frac{1}{k_1^2}\left(\left(\frac{k_3}{k_2}\right)'\right)^2 = constant. \tag{12}$$

Thus, using (11) and (12), we can characterize the spacelike $B_2$-slant helix $x(s)$ by the fact that

$$\left(\frac{k_3}{k_2}\right)^2 + \varepsilon_1 \frac{1}{k_1^2}\left(\left(\frac{k_3}{k_2}\right)'\right)^2 = constant \tag{13}$$

Conversely, if the condition (13) is satisfied for a regular spacelike curve we can always find a constant vector $U$ which makes a constant angle with the second binormal $B_2$ of the curve.

Consider the unit vector $U$ defined by

$$U = \left[\varepsilon_1 \frac{1}{k_1}\left(\frac{k_3}{k_2}\right)' T - \varepsilon_1 \frac{k_3}{k_2} N + B_2\right].$$

By taking account of the differentiation of (13), differentiation of $U$ gives that $\frac{dU}{ds} = 0$, this means that $U$ is a constant vector. So that, we can give the following theorem:

**Theorem 3.1.** *A unit speed spacelike curve $x: I \subset IR \to E_1^4$ with nonzero curvatures $k_1(s), k_2(s)$ and $k_3(s)$ is a spacelike $B_2$-slant helix if and only if the following condition is satisfied,*

$$\left(\frac{k_3}{k_2}\right)^2 + \varepsilon_1 \frac{1}{k_1^2}\left(\left(\frac{k_3}{k_2}\right)'\right)^2 = constant.$$

From Theorem 3.1 one can easily see that the constant function in Theorem 3.1, is independent of $\varepsilon_2$. So, we can give the following corollary.

**Corollary 3.2.** *The characterizations of the spacelike $B_2$-slant helix are independent of the Lorentzian causal character of the second binormal vector $B_2$. It is only related to the Lorentzian causal character of the unit principal normal vector $N$.*

Now, we give another characterization of spacelike $B_2$-slant helix in Minkowski 4-space.
Let assume that $x: I \subset IR \to E_1^4$ is a spacelike $B_2$-slant helix. Then, Theorem 3.1 is satisfied. By differentiating (13) with respect to $s$ we get

$$\left(\frac{k_3}{k_2}\right)\frac{d}{ds}\left(\frac{k_3}{k_2}\right) + \frac{\varepsilon_1}{k_1}\frac{d}{ds}\left(\frac{k_3}{k_2}\right)\frac{d}{ds}\left[\frac{1}{k_1}\frac{d}{ds}\left(\frac{k_3}{k_2}\right)\right] = 0, \qquad (14)$$

and hence

$$\frac{\varepsilon_1}{k_1}\left(\frac{k_3}{k_2}\right)' = -\frac{\left(\frac{k_3}{k_2}\right)\left(\frac{k_3}{k_2}\right)'}{\left[\frac{1}{k_1}\left(\frac{k_3}{k_2}\right)'\right]'}. \qquad (15)$$

If we write

$$f(s) = -\frac{\left(\frac{k_3}{k_2}\right)\left(\frac{k_3}{k_2}\right)'}{\left[\frac{1}{k_1}\left(\frac{k_3}{k_2}\right)'\right]'}, \qquad (16)$$

then

$$f(s)k_1 = \varepsilon_1\left(\frac{k_3}{k_2}\right)'. \qquad (17)$$

From (14) it can be written

$$\left[\frac{\varepsilon_1}{k_1}\left(\frac{k_3}{k_2}\right)'\right]' = -k_1\frac{k_3}{k_2}. \qquad (18)$$

By using (17) and (18) we have

$$\frac{d}{ds}f(s) = -k_1\frac{k_3}{k_2}. \qquad (19)$$

Conversely, let $f(s)k_1 = \varepsilon_1\left(\frac{k_3}{k_2}\right)'$ and $\frac{d}{ds}f(s) = -k_1\frac{k_3}{k_2}$. If we define a unit vector $U$ by

$$U = -f(s)T + \varepsilon_1\frac{k_3}{k_2}N - B_2 \qquad (20)$$

we have that $U$ and $\langle B_2, U \rangle$ are constants. So, we have the following theorem:

**Theorem 3.2.** *A unit speed spacelike curve $x: I \subset IR \to E_1^4$ with nonzero curvatures $k_1(s)$, $k_2(s)$ and $k_3(s)$ is a $B_2$-slant helix if and only if there exists a $C^2$-function $f$ such that*

$$fk_1 = \varepsilon_1 \frac{d}{ds}\left(\frac{k_3}{k_2}\right), \qquad \frac{d}{ds}f(s) = -k_1 \frac{k_3}{k_2}. \tag{21}$$

Now, we give the integral characterization of the spacelike $B_2$-slant helix.

Suppose that, the unit speed spacelike curve $x: I \subset IR \to E_1^4$ with nonzero curvatures $k_1(s), k_2(s)$ and $k_3(s)$ is a spacelike $B_2$-slant helix. Then theorem 3.2 is satisfied. Let us define $C^2$-function $\varphi$ and $C^1$-functions $m(s)$ and $n(s)$ by

$$\varphi = \varphi(s) = \int_0^s k_1(s) ds, \tag{22}$$

$$\begin{aligned} m(s) &= \frac{k_3}{k_2}\eta(\varphi) + f(s)\mu(\varphi), \\ n(s) &= \frac{k_3}{k_2}\mu(\varphi) - \varepsilon_1 f(s)\eta(\varphi), \end{aligned} \tag{23}$$

where $\eta(\varphi) = \cosh(\varphi)$, $\mu(\varphi) = \sinh(\varphi)$ if $\varepsilon_1 = -1$; and $\eta(\varphi) = \cos(\varphi)$, $\mu(\varphi) = \sin(\varphi)$ if $\varepsilon_1 = +1$. If we differentiate equations (23) with respect to $s$ and take account of (22) and (21) we find that $m' = 0$ and $n' = 0$. Therefore, $m(s) = C$, $n(s) = D$ are constants. Now substituting these in (23) and solving the resulting equations for $\frac{k_3}{k_2}$, we get

$$\frac{k_3}{k_2} = C\eta(\varphi) + D\mu(\varphi). \tag{24}$$

Conversely if (24) holds then from the equations in (23) we get
$$f = \varepsilon_1 (C\mu(\varphi) - D\eta(\varphi)),$$
which satisfies the conditions of Theorem 3.2. So, we have the following theorem:

**Theorem 3.3.** *A unit speed spacelike curve $x: I \subset IR \to E_1^4$ with nonzero curvatures $k_1(s), k_2(s)$ and $k_3(s)$ is a spacelike $B_2$-slant helix if and only if the following condition is satisfied*
$$\frac{k_3}{k_2} = C\eta(\varphi) + D\mu(\varphi),$$
*where $C$ and $D$ are constants, $\eta(\varphi) = \cosh(\varphi)$, $\mu(\varphi) = \sinh(\varphi)$ if $\varepsilon_1 = -1$; and $\eta(\varphi) = \cos(\varphi)$, $\mu(\varphi) = \sin(\varphi)$ if $\varepsilon_1 = +1$.*

**Conclusions**
In this paper, the spacelike $B_2$-slant helix is defined and the characterizations of the spacelike $B_2$-slant helix are given in Minkowski 4-space $E_1^4$. It is shown that a spacelike curve $x: I \subset IR \to E_1^4$ is a $B_2$-slant helix if an equation holds between the first, second and third curvatures of the curve. Furthermore, the integral characterization of the spacelike $B_2$-slant helix is given.